\documentclass{article}

\usepackage{amsmath,graphicx,amssymb}
\usepackage[utf8]{inputenc}
\usepackage{amssymb,amsmath,amscd,amsfonts,amsthm,bbm,xcolor,authblk}

\usepackage{caption}
\usepackage{subcaption}

\newcommand{\bs}{\boldsymbol}
\newcommand{\Sur}{{\mathcal S}}
\newcommand{\phat}{\hat{p}}
\newcommand{\shat}{\hat{\sigma}}
\newcommand{\mhat}{\hat{m}}

\newcommand{\jamal}[1]{ \color{red}{[Jamal:\ #1]}\ \color{black}}
\newcommand{\maxime}[1]{\color{blue}{[Max: #1]}\color{black}}

\newtheorem{theo}{Theorem}
\newtheorem{heur}{Heuristics}
\newtheorem{rem}{Remark}

\begin{document}
\title{Equilibrium and surviving species in a large Lotka-Volterra system of differential equations}

%
\author{Maxime Clenet$^{(1)}$, François Massol$^{(2)}$, Jamal Najim$^{(1)}$\thanks{Supported by CNRS Project 80 Prime | KARATE.}}
\affil{{\small (1)} CNRS and Université Gustave Eiffel, France\\ {\small (2)} CNRS, Université de Lille, INSERM, CHU, Institut Pasteur Lille}
%
%
\maketitle
\begin{abstract}
Lotka-Volterra (LV) equations play a key role in the mathematical modeling of various ecological, biological and chemical systems. When the number of species (or, depending on the viewpoint, chemical components) becomes large, basic but fundamental questions such as computing the number of surviving species still lack theoretical answers. In this paper, we consider a large system of LV equations where the interactions between the various species are a realization of a random matrix. We provide conditions to have a unique equilibrium and present a heuristics to compute the number of surviving species. This heuristics combines arguments from Random Matrix Theory, mathematical optimization (LCP), and standard extreme value theory. Numerical simulations, together with an empirical study where the strength of interactions evolves with time, illustrate the accuracy and scope of the results. 

\end{abstract}
%

%
\section{Introduction}
\label{sec:intro}

Since May's seminal work \cite{may1972will} and for the past decades, many theoretical studies addressed the issue of  the coexistence of species in ecosystems.

Introduced in the 1920s by Lotka \cite{lotka1925elements} and Volterra \cite{volterra1926fluctuations}, the Lotka-Volterra (LV) model is a well-known classic in theoretical ecology and mathematics. It represents a first step in our understanding of ecosystems through the variety of its dynamical behaviours (single or multiple equilibria, cycles, chaos),
its flexibility (many models can be approximated in the form of a LV model) and its mathematical calculability.

In this article, we consider large LV models with random parameters. Leveraging on the asymptotic understanding of large random matrices which naturally appear enables us to provide insights on equilibria and species coexistence for such models.   
\subsection*{Model and assumptions.} Large Lotka-Volterra systems of differential equations arise in various scientific fields such as biology, ecology, chemistry, etc. Although our results are generic in nature and not specific to a given field, we will rely on the ecological terminology in the sequel. 


A large system of Lotka-Volterra equations is a system of coupled ordinary differential equations (ODE) that write:
\begin{equation}\label{eq:LV}
\frac{dx_k(t)}{dt} = x_k(t)\, \left( r_k - \theta x_k(t) + \sum_{\ell\in [n]} B_{k\ell} x_{\ell}(t)\right)\ ,
\end{equation}
where $k\in [n]=\{1,\cdots, n\}$. 

Here, $n$ represents the number of species in a food web or community, the unknown vector $\boldsymbol{x}=(x_k)_{ k\in [n]}$ is the vector of abundances of the various species and evolves with time $t>0$ according to the dynamics \eqref{eq:LV}. 
Parameter $r_k$ represents the intrinsic growth rate of species $k$, $\theta$ is an intraspecific feedback coefficient (most often positive due to competition), and $B_{k\ell}$ is the per capita effect of species $\ell$ on species $k$. 

\begin{rem}
Notice that without interactions, i.e.
$B=(B_{k\ell})_{k,\ell\in [n]}=0\,,$
system \eqref{eq:LV} is simply a system of uncoupled logistic differential equations.
\end{rem}

We shall focus on the model where $r_k=\theta=1$:
\begin{equation}\label{eq:LV2}
\frac{dx_k}{dt} = x_k\, \left( 1 -  x_k + (B\boldsymbol{x})_k \right)\ ,\quad k\in [n]
\end{equation}
with matrix $B$ admitting the following representation:
\begin{equation*}
    B = \frac{A}{\alpha \sqrt{n}}+\frac{\mu}{n}\bs{1}\bs{1}^T\,,
\end{equation*}
where $A=(A_{ij})$ is a matrix with random standardized ($\mathbb{E}A_{ij}=0$ and $\mathrm{var}(A_{ij})=1$) independent and identically distributed (i.i.d.) entries with finite fourth moment, $\alpha>0$ is an extra parameter reflecting the interaction strength, and $\mu\in \mathbb{R}$ represents an arbitrary trend of the interactions. The $n\times 1$ vector $\bs{1}$ is a vector of ones.


\begin{rem}\label{rem:tao} Although matrix $B$ is a complex random object, a result by Tao \cite[Theorem 1.7]{tao_outliers_2013} fully describes its asymptotic spectrum: Assume that $|\mu|\, >\, 1/\alpha$, then for any fixed $\varepsilon>0$, almost surely eventually all the eigenvalues of $B$ but one are in the disk $\{z \in \mathbb{C}:|z| \leq 1/\alpha+\varepsilon \}$ while one extra eigenvalue takes the value $\mu+o(1)$.
\end{rem}

\begin{figure}
\centering
\begin{subfigure}[b]{0.48\textwidth}
  \centering
  \includegraphics[width=\textwidth]{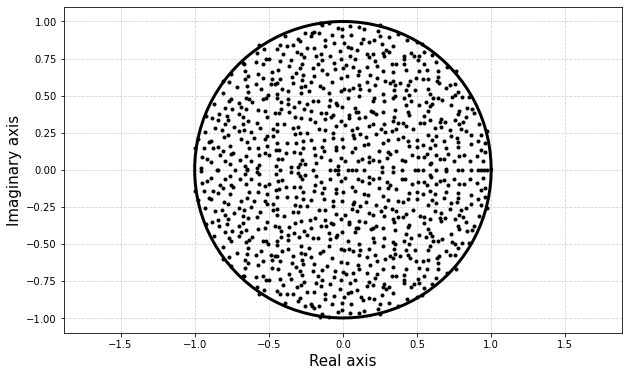}
  \caption{No outlier ($\mu=0$)}
  \label{subfig:circular}
\end{subfigure}
\hfill
\begin{subfigure}[b]{0.48\textwidth}
  \centering
  \includegraphics[width=\textwidth]{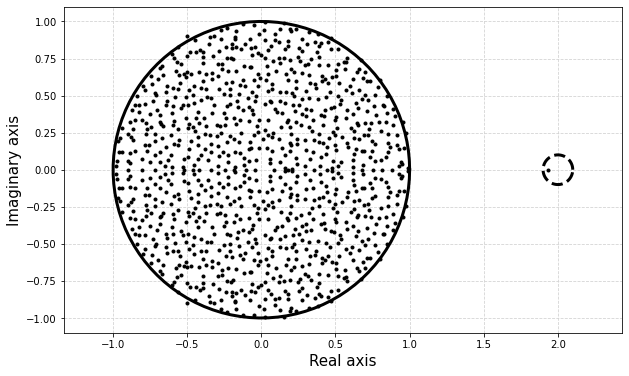}
  \caption{Presence of an outlier ($\mu=2$)}
  \label{subfig:outlier}
\end{subfigure}
\caption{Spectrum of non-Hermitian matrix B in the complex plan ($n=1000$, $\alpha=1$). In Fig. \ref{subfig:circular}, $\mu=0$ and the solid line circle represents the boundary of the circular law. In Fig. \ref{subfig:outlier}, $\mu=2$ and there is an eigenvalue in the small dashed circle centered at 2, as predicted by \cite[Th. 1.7]{tao_outliers_2013} - see also Remark \ref{rem:tao}.}
\label{fig:circ_law}
\end{figure}

Remark \ref{rem:tao} is illustrated in Fig. \ref{fig:circ_law}.


%

\subsection*{Presentation of the main results.}

\paragraph*{Unique equilibrium.} In the study of the behaviour of $\bs{x}(t)$ as $t\to \infty$ the existence of an equilibrium $\bs{x}^*$ to Eq. \eqref{eq:LV2} is an important prior to any stability property of $\bs{x}(t)$. By equilibrium, we mean the existence of a vector $\bs{x}^*=(x_k^*)$ satisfying $$
x^*_k(1-x^*_k+(B\bs{x}^*)_k)=0\quad \textrm{for}\quad k\in [n]\, .
$$
General results on LV systems state that $\bs{x}(t)>0$ (componentwise) as long as $\bs{x}(0)>0$ \cite{hofbauersigmund1998}. However, a possible equilibrium $\bs{x}^*$ will only verify $\bs{x}^*\ge 0$, i.e. some components $x_k^*$ will take the value zero.   

In Theorem \ref{th:unicite}, we provide sufficient conditions on the parameters $\alpha$ and $\mu$ to ensure the existence of a unique equilibrium. These conditions rely on the ``typical" behaviour of the random matrix $B$ in large dimension $n\to \infty$.

\paragraph*{Evaluating the number of surviving species.} Given a unique equilibrium $\bs{x}^*$, an important question is to describe the set of surviving/vanishing species. In this perspective, we introduce the set 
\begin{equation}\label{eq:def-S}
{\mathcal S}=\{ i\in [n],\  x^*_i>0\}
\end{equation}
of surviving species. In Section \ref{sec:heuristics}, we provide a heuristics to compute asymptotically the ratio $\frac {|{\mathcal S}|}n$ and understand via a system of equations the dependence between parameters $\alpha$ and $\mu$ and the number of surviving species. A complementary result addressing the elliptic random matrix model by means of theoretical physics methods can be found in \cite{bunin2017ecological} (dynamical cavity method) and in \cite{galla_dynamically_2018} (generating functional techniques).  

Notice that in \cite{bizeul2019positive}, Bizeul and Najim have studied a different normalization for $\alpha$ in the case $\mu = 0$, namely $\alpha \sim \sqrt{2\log(n)}$, to guarantee the survival of every species (feasibility of the equilibrium). Indeed, a consequence of Dougoud et al.'s results \cite{dougoud2018feasibility} is that some species will go to extinction if $\alpha>0$ is fixed (i.e. does not increase sufficiently with $n$). 

\paragraph*{An empirical study of LV systems with changing interaction strengths.}


Equipped with results on the existence of a unique equilibrium, one pending question is to understand what happens when the coefficient $\alpha$ varies for the same matrix $A$ and the same parameter $\mu$. In particular, when the value of $\alpha$ increases above a certain critical value, all species will coexist \cite{bizeul2019positive}; conversely, for sufficiently low values of $\alpha$, the existence of a feasible equilibrium is not warranted anymore, however a unique and stable equilibrium may exist. How species equilibrium abundances change between these two states and how $|\mathcal S|$ varies will be the focus of Section \ref{sec:ecological}.

\subsection*{Notations} Denote by $\rho(C)$ the spectral radius of matrix $C$, by $\|C \|$ the spectral norm of matrix $C$, and by $\| \bs{u}\|$ the euclidean norm of vector $\bs{u}$. We represent by $\delta_x$ the Dirac measure at $x$: 
$$
\delta_x(E)=\begin{cases} 1& \textrm{if}\ x\in E\\
0&\textrm{else}
\end{cases}\, .
$$
We denote by $\xrightarrow[]{a.s.}$ the almost sure convergence of random quantities and by $\xrightarrow[]{weak}$ the weak convergence of measures. Given a set $S$, we denote by $|S|$ its cardinality.

\section{Equilibrium and stability results}\label{sec:equilibrium}

\paragraph*{A primer on Random Matrix Theory.} We first recall some results on Random Matrix Theory (RMT), which provides a number of valuable insights to understand the asymptotic behaviour of $A$. We begin by the almost sure (a.s.) convergence of the spectral radius and the spectral norm:
$$
\rho(A/\sqrt{n})\xrightarrow [n\to\infty]{a.s.}1\quad \text{and} \quad \|A/\sqrt{n}\|\xrightarrow [n\to\infty]{a.s.}2\, .
$$ 
We also have the a.s. weak convergence of the spectral measure of $A/\sqrt{n}$ to the circular law (see for instance \cite{2012-bordenave-chafai-circular}):
$$
(a.s.)\quad \frac 1n \sum_{k\in [n]} \delta_{\lambda_k(A/\sqrt{n})} \ \xrightarrow[n\to\infty]{weak}\  \frac{\boldsymbol{1}_{\{x^2+y^2\le 1\}}}{\pi}\, dx\, dy\, ,
$$
where $(\lambda_k(A/\sqrt{n}); k\in [n])$ is the spectrum of $A/\sqrt{n}$. This convergence is illustrated in Fig. \ref{subfig:circular}.

The description of the spectral norm of the deterministic part of matrix $B$ is more straightforward:
$$
\left\| \mu \, \frac{\bs{1} \bs{1}^T}n\right\| =|\mu|\, .
$$
Notice that both the random and deterministic parts of matrix $B$ do not vanish asymptotically and thus have a macroscopic effect on the dynamics of system \eqref{eq:LV}, as recalled in Remark \ref{rem:tao} where the asymptotic spectrum of $B$ is described. 

\paragraph*{The non-invadability condition.}
A key element to understand the dynamics of the LV system \eqref{eq:LV} is the existence of an equilibrium $\boldsymbol{x}^*=(x^*_k)_{k\in[n]}$ such that 
\begin{equation}\label{eq:equilibrium}
\begin{cases}
 & x^*_k\left( 1 -  x^*_k + (B\boldsymbol{x^*})_k \right) =0\,,\quad \forall k\in [n]\, ,\\ 
 & x^*_k \geq 0.
\end{cases}
\end{equation}
and the study of its stability, that is the convergence of a solution $\bs{x}$ to the equilibrium $\bs{x}^*$: $\bs{x}(t)\xrightarrow[t\to\infty]{} \bs{x}^*$ if $\bs{x}(0)$ is sufficiently close to $\bs{x}^*$.

It is well known that for LV equations, the fact that $\bs{x}(0)>0$ (componentwise) implies that $\bs{x}(t)>0$ for every $t>0$, but one can have some components $x_k(t)$ of $\bs{x}(t)$ vanishing to zero. As a consequence, we will only consider non negative equilibria $\bs{x}^*\ge 0$ with possibly vanishing components.

Notice that the situation substantially differs whether $\bs{x}^*>0$ or $\bs{x}^*$ has vanishing component. In the former case, the equilibrium set of equations becomes a linear equation: 
$$
\bs{x}^* = \bs{1}+B  \bs{x}^*\, .
$$ 
In the latter case, the equilibrium equations are no longer linear. 

In the centered case $\mu=0$, the existence of a positive solution has been studied in \cite{bizeul2019positive} and requires $\alpha\gg \sqrt{2\log(n)}$ (while we consider $\alpha$ fixed here).

A naive and systematic way to solve \eqref{eq:equilibrium} is to choose a priori a subset ${\mathcal I}\subset [n]$, to set the corresponding components $\bs{x}_{\mathcal I}=(x_i^*)_{i\in {\mathcal I}}$ to zero, and to solve the remaining linear system:
$$
\bs{x}_{\mathcal{I}^c} = \bs{1}_{|\mathcal{I}^c|} + B_{\mathcal{I}^c}\bs{x}_{\mathcal{I}^c}\ .
$$
If there exists $\bs{x}_{\mathcal{I}^c}\ge 0$ that solves the previous equation, then $\bs{x}=\begin{pmatrix}\bs{x}_{\mathcal I}\\
\bs{x}_{\mathcal{I}^c}
\end{pmatrix}$ satisfies \eqref{eq:equilibrium} and is a potential equilibrium. The number of subcases ${\mathcal I} \subset [n]$ is $2^n$ and in particular grows exponentially as $n\to\infty$.


In order to decrease the number of potential solutions to \eqref{eq:equilibrium}, we first notice that 
relying on standard properties of dynamical systems, see for instance \cite[Theorem 3.2.5]{takeuchi1996global}, a necessary condition for the equilibrium $\bs{x}^*$ to be stable is that 
\begin{equation}
\label{non_inva_1}
1-x^*_k+(B  \bs{x}^*)_k\le 0\, .
\end{equation} 
The condition (\ref{non_inva_1}) is better known in ecology as the non-invadability condition \cite{law_permanence_1996}. In reference to the ODE (\ref{eq:LV2}), the requirement for a given species $k \in [n]$ to be non-invasive is equivalent to:
\begin{equation}
    \label{non_inva_2}
    \left( \frac{1}{x_k} \frac{dx_k}{dt}\right)_{x_k \rightarrow 0^+}\le 0\, .
    \end{equation}
The main interpretation is as follows: if one adds species $k$ with a very low abundance in the system, it will not be able to invade the system as a result of condition (\ref{non_inva_2}). 

As a consequence, we will now focus on the following set of conditions:
\begin{equation}\label{eq:equilibrium-NI}
\left\{
\begin{array}{cccl}
 x^*_k\left( 1 -  x^*_k + (B\boldsymbol{x^*})_k \right) &=&0&\textrm{for}\ k\in [n]\, ,\\ 
  1 -  x^*_k + (B\boldsymbol{x^*})_k  &\le& 0&\textrm{for}\ k\in [n]\, ,\\
 \bs{x}^* &\geq& 0 & \textrm{componentwise}\, .\\
 \end{array}
\right.
\end{equation}

This casts the problem of finding a non negative equilibrium into the class of Linear Complementarity Problems (LCP), which we describe hereafter.

\paragraph*{Linear Complementarity Problem (LCP).} LCP is a class of problems from mathematical optimization which in particular encompasses linear and quadratic programs; standard references are \cite{murty1988linear,cottle2009linear}.  Given a $n\times n$ matrix $M$ and a $n\times 1$ vector $\bs{q}$, the associated LCP denoted by $LCP(M,\bs{q})$ consists in finding two $n\times 1$ vectors $\bs{z},\bs{w}$ satisfying the following set of constraints:
\begin{equation}\label{eq:solution-LCP}
\left\{
\begin{array}{lcl}
\bs{z}&\ge& 0\, , \\
\bs{w}=M\bs{z} +\bs{q}&\ge& 0\, , \\
\bs{w}^T \bs{z}&=&0\quad  \Leftrightarrow\quad w_kz_k=0\quad\ \textrm{for all}\ k\in [n]\,.
\end{array}
\right.
\end{equation}
Since $\bs{w}$ can be inferred from $\bs{z}$, we denote $\bs{z}\in LCP(M,\bs{q})$ if $(\bs{w},\bs{z})$ is a solution of \eqref{eq:solution-LCP}.

A theorem by Murty \cite{murty1972number} states that the $LCP(M,\bs{q})$ has a unique solution $(\bs{w},\bs{z})$ iff $M$ is a $P$-matrix, that is:
$$
\det(M_{\mathcal I}) >0\ , \quad \forall\, {\mathcal I}\subset [n]\ ,\quad M_{\mathcal I}=(M_{k\ell})_{k,\ell\in {\mathcal I}}\, .
$$
In view of \eqref{eq:equilibrium-NI}, we look for $\bs{x}^*\in LCP(I-B, -\bs{1})$.

\paragraph*{The equilibrium $\bs{x}^*$ and its stability.} For a generic LV system
\begin{equation}\label{eq:generic-LV}
\frac{d\, y_k(t)}{dt} =y_k(r_k +(C\bs{y})_k)\, ,\quad k\in [n]\, ,
\end{equation}
Takeuchi and Adachi (see for instance \cite[Th. 3.2.1]{takeuchi1996global}) provide a criterion for the existence of a unique equilibrium $\bs{y}^*$ and the global stability of the LV system.
\begin{theo}[Takeuchi and Adachi \cite{takeuchi1980existence}] \label{th:takeuchi} If there exists a positive diagonal matrix $\Delta$ such that $\Delta C + C^T \Delta$ is negative definite, then $LCP(-C,\bs{r})$ admits a unique solution. In particular, for every $\bs{r}\in \mathbb{R}^n$, there is a unique equilibrium $\bs{y}^*$ to \eqref{eq:generic-LV}, which is globally stable in the sense that for every $\bs{y}_0>0$, the solution to \eqref{eq:generic-LV} which starts at $\bs{y}(0)=\bs{y}_0$ satisfies
$$
\bs{y}(t)\xrightarrow[t\to\infty]{} \bs{y}^*\, .
$$
\end{theo} 
Combining this result (setting $C=-(I-B)$) with results from RMT, we can guarantee the existence of a globally stable equilibrium $\bs{x}^*$ of \eqref{eq:LV} for a wide range of the set $(\alpha,\mu)$. Denote by 
\begin{equation}\label{eq:def-A}
{\mathcal A}=\left \{ (a,m)\in \mathbb{R}^*_+ \times \mathbb{R}: \ a > \sqrt{2},\ m< \frac{1}{2}+\frac 12 \sqrt{1-\frac{2}{a^2}}  \right \}
\end{equation}
the set of {\it admissible parameters}.
\begin{theo}\label{th:unicite} Let $(\alpha,\mu) \in {\mathcal A}$, then a.s. matrix 
$
(I-B)+(I-B)^T
$
is eventually positive definite: with probability one, for a given realization $\omega$, there exists $N(\omega)$ such that for $n\ge N(\omega)$, $(I-B^\omega)+(I-B^\omega)^T$ is positive definite. In particular, there exists a unique (random) globally stable equilibrium $\bs{x}^*\in LCP(I-B^\omega, -\bs{1})$ to \eqref{eq:equilibrium-NI}.
\end{theo}

\begin{figure}[h!]
  \centering
  \centerline{\includegraphics[width=8.5cm]{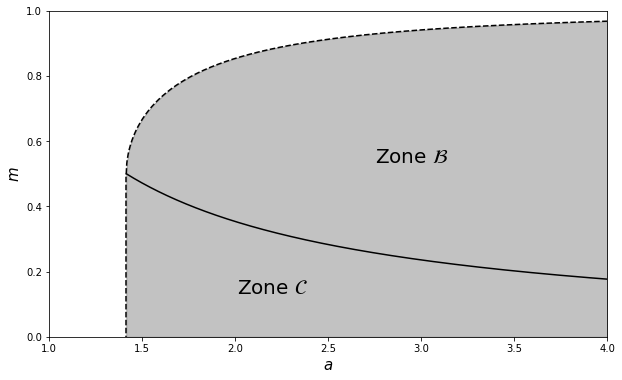}}
   \caption{The shaded area represents the set $\mathcal{A}$ given by \eqref{eq:def-A} yielding the existence of a unique (random) globally stable equilibrium $\bs{x}^*$. Area ${\mathcal A}$ is divided into two zones ${\mathcal B}$ and ${\mathcal C}$. Both zones correspond to parameters ($\alpha,\mu$) for which matrix $2I -(B+B^T)$ is definite positive, as stated in Theorem \ref{th:unicite}. In zone ${\mathcal B}$, $\lambda_{\max}(B+B^T)$ corresponds to a spiked eigenvalue ($\mu$ above the critical threshold $(\alpha\sqrt{2})^{-1}$). In zone ${\mathcal C}$, $\lambda_{\max}(B+B^T)$ corresponds to the right edge of the semi-circle law. Notice that zone ${\mathcal C}$ extends to negative values along the $y$-axis.}  
\label{fig:admissible}
\end{figure}

\begin{figure}
\centering
\begin{subfigure}[b]{.48\textwidth}
    \includegraphics[width=\textwidth]{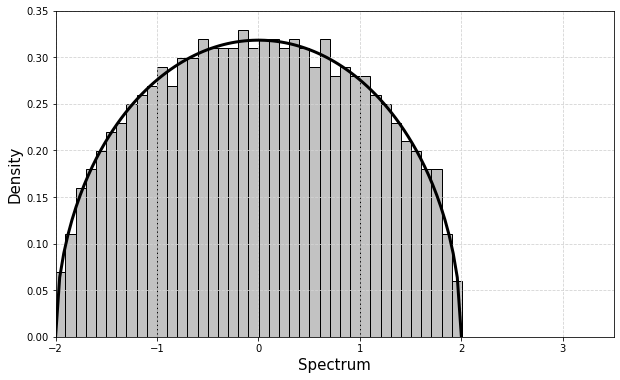}
    \caption{No outlier if $\mu\le (\alpha\sqrt{2})^{-1}$. }
    \label{subfig:SC}
  \end{subfigure}%
  \hfill   
  \begin{subfigure}[b]{0.48\textwidth}
    \includegraphics[width=\textwidth]{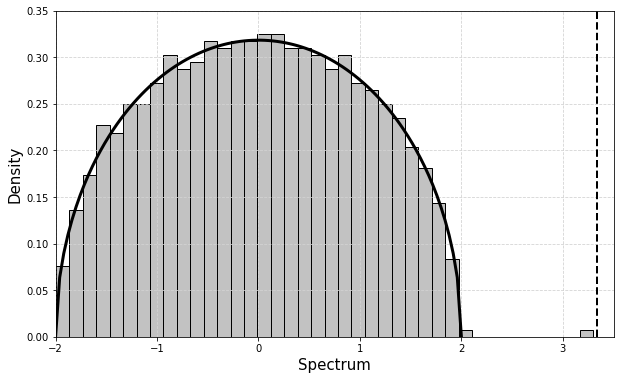}
    \caption{Outlier if $\mu> (\alpha\sqrt{2})^{-1}$.}
    \label{subfig:SC-outlier}
  \end{subfigure}%
\caption{Spectrum (histogram) of the Hermitian random matrix $B+B^T$ ($n=1000$, $\alpha=\sqrt{2}$). The solid line represents the semi-circular law. In Fig. \ref{subfig:SC}, $\mu=0$ and there is no oulier. In Fig. \ref{subfig:SC-outlier}, $\mu=1.5$ and one can notice the presence of an eigenvalue outside the bulk of the circular law. The dashed line indicates its theoretical value.}
\label{fig:semi-circle}
\end{figure}

\begin{proof}
We have
$$
I- B +I -B^T = 2I - (B+B^T) = 2I - \left(\frac{A+A^T}{\alpha \sqrt{n}}+\frac{2\mu}{n} \bs{1}\bs{1}^*\right)\ .
$$
Notice that $2I- (B+B^T)$ is positive definite iff the top eigenvalue of $B+B^T$ is lower than 2:
\begin{equation}
\label{cond_eig}
\lambda_{\max}(B+B^T) < 2 
\end{equation}
We first focus on the random part $(A+A^T)/\alpha$ which is a symmetric matrix with independent ${\mathcal N}(0,2/\alpha^2)$ entries above the diagonal (note that the distribution of the diagonal entries is different from the off-diagonal entries, with no asymptotic effect). In this case, it is well known that the largest eigenvalue of the normalized matrix (or equivalently its spectral norm since the matrix is symmetric) a.s. converges to the right edge of the support of the semi-circle law (see \cite[Th. 5.2]{book-bai-silverstein}):
\begin{equation}
\label{centered_cond}
\lambda_{\max}\left( \frac{A+A^T}{\alpha \sqrt{n}}
\right)
\xrightarrow[n\to\infty]{a.s.} \frac{2 \sqrt{2}}{\alpha}\, .
\end{equation}
In the centered case ($\mu = 0$), condition (\ref{cond_eig}) occurs if $\alpha > \sqrt{2}$.

We now consider the general case where $\mu \neq 0$. Notice that the rank-one perturbation matrix $P = \frac{2\mu}{n} \bs{1}\bs{1}^*$ admits a unique non zero eigenvalue $2\mu$. Denote by $\check{A}=\frac{A+A^T}{\alpha\sqrt{n}}$. We are concerned with the top eigenvalue of the symmetric matrix $\check{A}+P$. Based on a result by Capitaine et al. \cite[Th. 2.1]{capitaine_largest_2009}, we have:
\begin{equation*}
    \lambda_{\max}(\check{A}+P) \xrightarrow[n \rightarrow \infty]{a.s}\left\{ 
 \begin{array}{ll}
  2\mu+\frac{1}{\alpha^2  \mu} &\text{if } \mu> \frac{1}{\sqrt{2} \alpha}\, ,  \\ 
  \frac{2 \sqrt{2}}{\alpha} &\text{else.}
\end{array}\right.
\end{equation*}
This result is illustrated in Figure \ref{fig:semi-circle}.

Assume first that $\mu\le \frac 1{\alpha\sqrt{2}}$ (corresponding to zone ${\mathcal C}$ in Fig. \ref{fig:admissible}), then $\lambda_{\max}( \check{A}+P) \xrightarrow[n\to\infty]{a.s.} \frac{2\sqrt{2}}{\alpha}$, which is strictly lower than 2 (cf. condition \eqref{cond_eig}) if $\alpha>\sqrt{2}$. Hence $\lambda_{\max}(\check{A} +P)$ is eventually strictly lower than 2 under this condition.

Assume now that $\mu>\frac 1{\alpha\sqrt{2}}$ (corresponding to zone ${\mathcal B}$ in Fig. \ref{fig:admissible}), then 
$$
\lambda_{\max}( \check{A}+P) \xrightarrow[n\to\infty]{a.s.} 2\mu+\frac{1}{\alpha^2  \mu}\, .
$$
We are interested in the conditions for which $2\mu +\frac 1{\alpha^2 \mu}<2$ or equivalently 
\begin{equation}\label{eq:condition-poly}
    2\alpha^2 \mu^2 - 2\alpha^2 \mu +1<0\, .
\end{equation}

An elementary study of the polynomial $\xi(X)= 2\alpha^2 X^2 - 2\alpha^2 X +1$ yields that $\xi$'s discriminant is positive if $\alpha>\sqrt{2}$, $$\xi(\mu^{\pm})=0\quad \Leftrightarrow\quad  \mu^{\pm}= \frac 12 \pm \frac 12 \sqrt{1- \frac{2}{\alpha^2}}\ ,
$$ 
and $\xi\left( \frac 1{\alpha\sqrt{2}}\right) <0$, so that $\frac 1{\alpha\sqrt{2}}\in (\mu^-, \mu^+)$. In particular condition \eqref{eq:condition-poly} is fulfilled if 
$$
\mu \in \left( \frac 1{\alpha\sqrt{2}}\ ,\ \frac{1}{2}+\frac 12 \sqrt{1-\frac{2}{\alpha^2}}\right)\ .
$$
Under this condition, \eqref{eq:condition-poly} is fulfilled and a.s. $\limsup_{n\to\infty}\lambda_{\max}(\check{A}+P) < 2$, which completes the proof: we can then rely on Theorem \ref{th:takeuchi} to conclude.
\end{proof}

%

\section{A heuristic approach to the proportion and distribution of the surviving species}
\label{sec:heuristics}
\subsection{Proportion of surviving species}

In Section \ref{sec:equilibrium}, we have presented conditions on parameters $\alpha, \mu$ for the existence of a globally stable equilibrium 
$\bs{x}^*$ to \eqref{eq:LV} under the non-invadability condition. As $\bs{x}^*$ depends on the realization of matrix $B$, it is a random vector. Moreover since $\alpha>0$ is fixed and does not depend on $n$, the equilibrium $\bs{x}^*$ will feature vanishing components (see the original argument in \cite{dougoud2018feasibility} and the discussion in \cite{bizeul2019positive}).
In an ecological context, we shall refer to these non-vanishing components as the surviving species, the vanishing components corresponding to the species going to extinction with $x_k^*=0$ and
$
x_k(t)\xrightarrow[t\to\infty]{} 0\, .
$

In this section, we assume that the $A_{ij}$'s are ${\mathcal N}(0,1)$-distributed and describe the proportion of non-vanishing components of the equilibrium $\bs{x}^*$; we also describe the distribution of the surviving species $x_i^*>0$ which turns out to be a truncated Gaussian. 

\begin{rem}
The Gaussianity assumption facilitates the explanation of the heuristics but does not seem necessary for the result to hold. In Fig. \ref{subfig:unif}, the entries are not considered Gaussian but the distribution of the surviving species still matches the truncated Gaussian.
\end{rem}

Given the random equilibrium $\bs{x}^*$, recall the definition of ${\mathcal S}$ in \eqref{eq:def-S}. We introduce the following quantities:
$$
\phat=\frac {|\Sur|}{n}\ ,\qquad \mhat = \frac 1{|\Sur|} \sum_{i\in [n]} x_i^*\ ,\qquad \shat^2 =\frac 1{|\Sur|} \sum_{i\in [n]} (x_i^*)^2\,.
$$
Notice that in the definitions of $\mhat$ and $\shat^*$ we can replace $\sum_{i\in [n]}$ by $\sum_{i\in {\mathcal S}}$.

Denote by $Z\sim{\mathcal N}(0,1)$ a standard Gaussian random variable and by $\Phi$ the cumulative Gaussian distribution function:
$$
\Phi(x)=\int_{-\infty}^x \frac{e^{-\frac{u^2}2}}{\sqrt{2\pi}}\, du\, .
$$
Recall the definition of the set ${\mathcal A}$ in \eqref{eq:def-A}.
\begin{heur}\label{heur:proportion}  Let $(\alpha,\mu) \in {\mathcal A}$. The following system of three equations and three unknowns $(p ,m,\sigma)$
\begin{eqnarray}\label{eq:heur1}
\sigma\sqrt{p} \Phi^{-1}(1-p) +\alpha(1+\mu\, p\, m) &=& 0\,,\\
\label{eq:heur2}
1 + \mu\, p\, m +\frac{\sigma\sqrt{p}}{\alpha} \mathbb{E}(Z\mid Z>-\delta) &=& m\, , \\
\nonumber
(1+ \mu\, p\, m)^2+(1+ \mu\, p\, m)\frac{2\sigma \sqrt{p}}{\alpha}\mathbb{E}(Z\mid Z>-\delta) \quad \qquad &\phantom{=}&\\
\label{eq:heur3}
+\quad  \frac{\sigma^2 p}{\alpha^2} \mathbb{E} (Z^2\mid Z> -\delta)   &=& \sigma^2
\end{eqnarray}
where \begin{equation}\label{eq:def-delta}
\delta = \delta(p,m,\sigma) = \frac{\alpha}{\sigma\sqrt{p}}(1+ \mu \pi m)\ ,
\end{equation}
admits a unique solution $(p^*,m^*,\sigma^*)$
 and
 $$
 \phat \xrightarrow[n\to\infty]{a.s.} p^*\ ,\qquad \mhat\xrightarrow[n\to\infty]{a.s.} m^*\,\qquad \text{and}\qquad \shat\xrightarrow[n\to\infty]{a.s.} \sigma^*\, .
$$ 
Associated to this solution $(p^*,m^*,\sigma^*)$ is $\delta^*=\delta(p^*,m^*,\sigma^*)$.
\end{heur}

\begin{figure}[h!]
\centering
\begin{subfigure}{0.48\textwidth}
  \includegraphics[width=\textwidth]{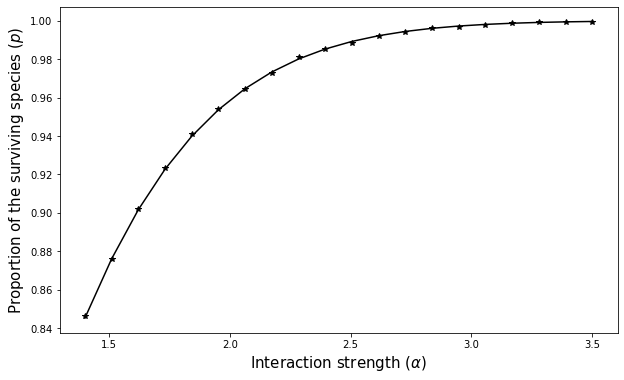}
  \label{}
\end{subfigure}
\hfill
\begin{subfigure}{0.48\textwidth}
  \includegraphics[width=\textwidth]{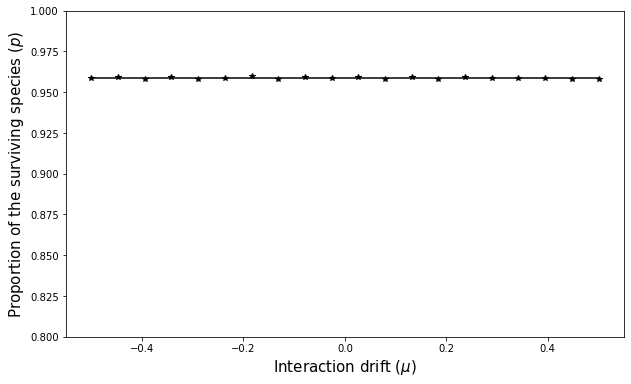}
  \label{}
\end{subfigure}
\hfill
\begin{subfigure}{0.48\textwidth}
  \includegraphics[width=\textwidth]{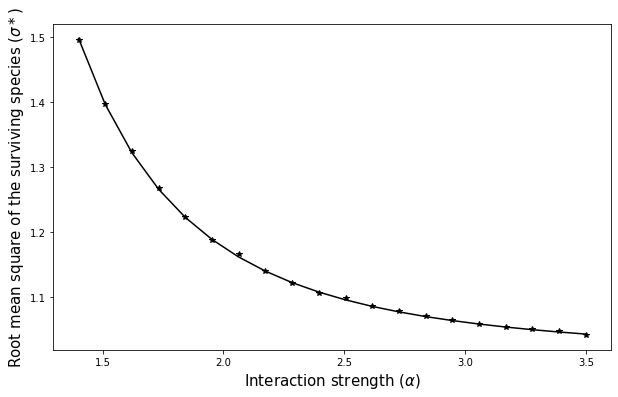}
  \label{}
\end{subfigure}
\hfill
\begin{subfigure}{0.48\textwidth}
  \includegraphics[width=\textwidth]{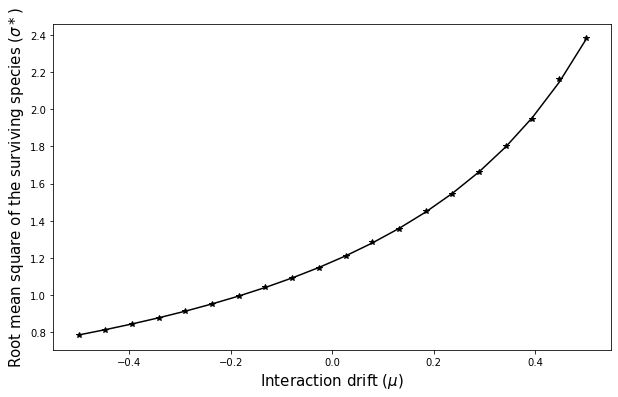}
  \label{}
\end{subfigure}
\hfill
\begin{subfigure}{0.48\textwidth}
  \includegraphics[width=\textwidth]{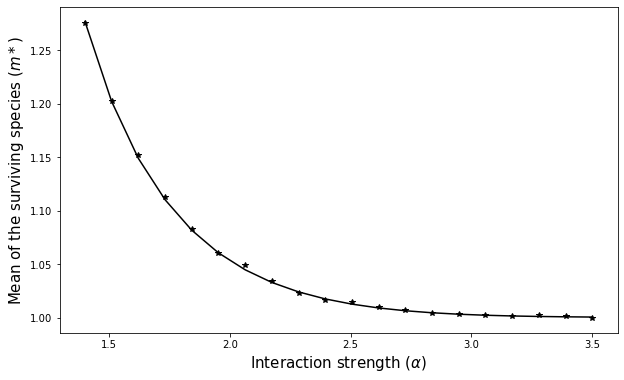}
  \caption{Parameters $(p^*,\sigma^*,m^*)$ versus $\alpha$.}
  \label{subfig:vsalpha}
\end{subfigure}
\hfill
\begin{subfigure}{0.48\textwidth}
  \includegraphics[width=\textwidth]{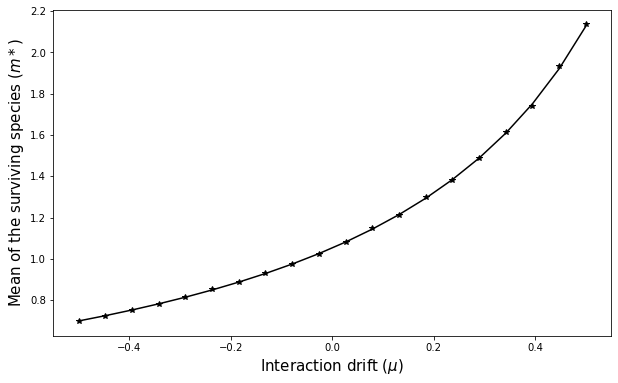}
  \caption{Parameters $(p^*,\sigma^*,m^*)$ versus $\mu$.}
  \label{subfig:vsmu}
\end{subfigure}
\caption{The plots represent a comparison between the theoretical solutions $(p^*,\sigma^*,m^*)$ of \eqref{eq:heur1}-\eqref{eq:heur3}  and their empirical Monte Carlo counterpart (the star marker) as functions of the interaction strength $\alpha$ (left) and the interaction drift $\mu$ (right). Matrix B has size $n=500$ and the number of Monte Carlo experiments is 200. In Column \eqref{subfig:vsalpha}, $\mu=0$ and $\alpha > \sqrt{2}$ on the $x$-axis (which guarantees a unique and stable equilibrium $\bs{x}^*$). When interaction $\alpha^{-1}$ increases, the number of surviving species $p^*$ decrease but their variance $\sigma^*$ and mean $m^*$ increase.
In Column \eqref{subfig:vsmu}, $\alpha=2$ and $\mu\in (-0.5,0.5)$ on the $x$-axis. The interaction drift appears to have no impact on the proportion $p^*$ of surviving species, whereas it influences their variances and means.}
\label{fig:test}
\end{figure}

\begin{figure}[t!]
\begin{minipage}[b]{0.98\linewidth}
  \centering
  \centerline{\includegraphics[width=8.5cm]{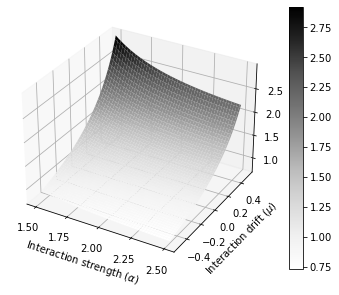}}
   \caption{The 3D plot represents $\sigma^*=\sigma^*(\alpha,\mu)$, solution of the system \eqref{eq:heur1}-\eqref{eq:heur3}. In contrast to the proportion of surviving species $p^*$, we observe that $\mu$ has a major influence over $\sigma^*$. The graph for the theoretical value of $m^*$ has approximately the same behavior with respect to $\mu$ and $\alpha$. 
   }
\label{fig:sensitivity}
\end{minipage}
\end{figure}

There is a strong matching between the parameters obtained by solving \eqref{eq:heur1}-\eqref{eq:heur3} and their empirical counterparts obtained by Monte-Carlo simulations. This is illustrated in Fig. \ref{fig:test}. In Fig. \ref{fig:sensitivity}, we illustrate the sensitivity of $\sigma^*$ to the parameters $(\alpha,\mu)$. 

\begin{rem}
The heuristics above substantially simplifies in the centered model case, where $\mu=0$ and $B=\frac{A}{\alpha\sqrt{n}}$. Following \eqref{centered_cond}, assume that $\alpha>\sqrt{2}$. Then the system with two unknowns $(p ,\sigma)$
$$
\left\{ 
\begin{array}{lcl}
\sigma\sqrt{p} \Phi^{-1}(1-p) +\alpha\phantom{\bigg|} &=& 0\\
1+\frac{2\sigma \sqrt{p}}{\alpha}\mathbb{E}(Z\mid Z>-\delta) +  \frac{\sigma^2 p}{\alpha^2} \mathbb{E} (Z^2\mid Z> -\delta)   &=& \sigma^2
\end{array}\right. \, \textrm{where}\quad \delta  = \frac{\alpha}{\sigma\sqrt{p}}
$$
admits a unique solution $(p^*,\sigma^*)$. Moreover, $\phat \xrightarrow[n\to\infty]{a.s.} p^*$ and $\shat\xrightarrow[n\to\infty]{a.s.} \sigma^*$.
\end{rem}

\subsection{Distribution of surviving species}

In the previous section, the proportion $\phat$ of the surviving species, their mean $\mhat$ and second moment $\shat^2$ have been described as empirical counterparts of the solutions $p^*, m^*, (\sigma^*)^2$ of a system of equations. 
While establishing this system of equations, we will provide the following representation (see \eqref{eq:representation-x_k}) of the abundance $x^*_k$ of a surviving species:
$$
x_k^* \ =\ 1+\mu\,p^* m^* +\frac{\sigma^*\sqrt{p^*}}{\alpha} Z_k\,,
$$
where $Z_k\sim{\mathcal N}(0,1)$ and $Z_k>-\delta^*=-\delta(p^*,m^*,\sigma^*)$, $\delta$ being defined in \eqref{eq:def-delta}. We take here advantage of this representation 
to characterize $x_k$'s distribution, which turns out to be a truncated Gaussian.

\begin{heur}\label{heur:distribution} Let $(\alpha,\mu) \in {\mathcal A}$, $\bs{x}^*$ the solution of \eqref{eq:equilibrium-NI} and let $(p^*,m^*,\sigma^*)$ the solution of the system \eqref{eq:heur1}-\eqref{eq:heur3}. Recall the definition \eqref{eq:def-delta} of $\delta$ and denote by $\delta^*=\delta(m^*,p^*,\sigma^*)$. Let $x^*_k>0$ a positive component of $\bs{x}^*$, then:
$$
\mathcal{L}(x_k^*) \xrightarrow[n\to\infty]{} \mathcal{L}\left(1+\mu p^* m^* +\frac{\sigma^*\sqrt{p^*}}{\alpha} Z \quad \bigg|\quad Z > -\delta^* \right)\ ,
$$
where $Z\sim {\mathcal N}(0,1)$. Otherwise stated, asymptotically $x^*_k$ admits the following density
$$
 f(y) = \frac{\bs{1}_{\{y>0\}}}{\Phi(\delta^*)}\frac{\alpha}{\sigma^*\sqrt{2\pi\,p^*}}\,
 \exp\left\{ - -\frac 12\left(\frac{\alpha}{\sigma^* \sqrt{p^*}}y-\delta^*\right)^2\right\}\ .
$$
\end{heur}

The heuristics simply follows from the fact that if $x_k^*$ is a surviving species then
$$
x_k^* = 1+\mu\,p^* m^* +\frac{\sigma^*\sqrt{p^*}}{\alpha} Z_k
$$
conditionally on the fact that the right hand side of the equation is positive, that is $Z_k>-\delta^*$. A simple change of variable yields the density - details are provided in Appendix \ref{app:proof-heuristics}.

Fig. \ref{fig:histo} illustrates the matching between the theoretical distribution obtained in Heuristics \ref{heur:distribution} and a histogram obtained by Monte-Carlo simulations. It also illustrates the validity of the heuristics beyond the gaussiannity assumption of the entries.

\begin{figure}
\centering
    \begin{subfigure}[b]{0.48\textwidth}
    \centering
    \includegraphics[width=\textwidth]{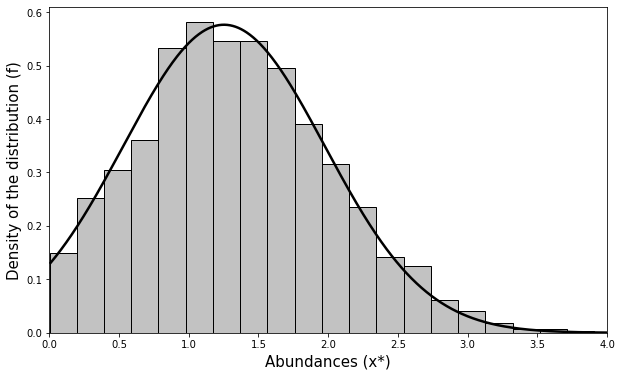}
    \caption{Gaussian entries.}\label{subfig:gauss}
    \end{subfigure}
    \hfill
    \begin{subfigure}[b]{0.48\textwidth}
    \centering
    \includegraphics[width=\textwidth]{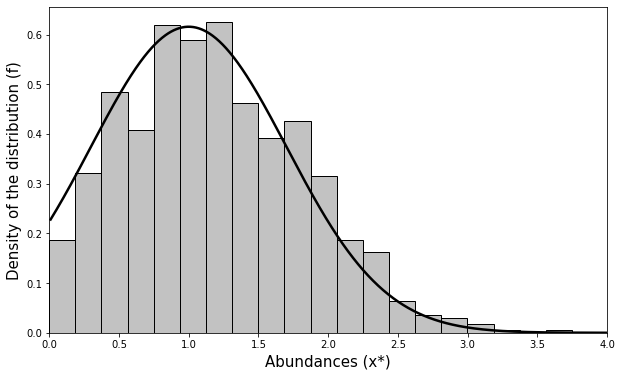}
    \caption{Uniform entries.}
    \label{subfig:unif}
    \end{subfigure}
  \caption{Distribution of surviving species. The $x$-axis represents the value of the abundances and the histogram is built upon the positive components of equilibrium $\bs{x}^*$. The solid line represents the theoretical distribution for parameters $(\alpha,\mu)$ as given by Heuristics \ref{heur:distribution}. In Fig. \eqref{subfig:gauss}, the entries are Gaussian ${\mathcal N}(0,1)$ and the parameters are set to $(n=2000,\alpha = 2,\mu = 0.2)$. In Fig. \eqref{subfig:unif}, the entries are uniform $\mathcal{U}(-\sqrt{3},\sqrt{3})$ with variance 1 and the parameters are set to $(n=2000,\alpha = \sqrt{3},\mu = 0)$. Notice in particular that the theoretical distribution matches with non-Gaussian entries.}
  \label{fig:histo}
\end{figure}

\subsection{Construction of the heuristics}
We first discuss Heuristics \ref{heur:proportion} and establish Equations \eqref{eq:heur1}, \eqref{eq:heur2} and \eqref{eq:heur3}.
\paragraph*{Equation \eqref{eq:heur1}.} We first recall a result on order statistics of a Gaussian sample. Consider a family $(Z_k)_{k\in [n]}$ of i.i.d. random variables ${\mathcal N}(0,1)$ and the associated order statistics
$$
Z_1^*\ \le\ Z_2^*\ \le \cdots \le\ Z_n^*\, .
$$  
Consider an index $\lfloor n\alpha \rfloor\in [n]$ where $\alpha\in (0,1)$ is fixed, then the typical location of $Z^*_{ \lfloor n\alpha \rfloor}$ is $\Phi^{-1}(\alpha)$:
\begin{equation}\label{eq:typical}
Z^*_{ \lfloor n\alpha \rfloor}\simeq \Phi^{-1}(\alpha)\quad \text{as}\quad n\to\infty\, ,
\end{equation}
see for instance \cite{smirnov1949limit,balkema1978limit}. 

Let $\bs{x}^*$ be the equilibrium of \eqref{eq:LV} and consider the random variable 
$$\check{Z}_k= \sum_{i\in {\mathcal S}} B_{ki} x_i^* = (B\bs{x}^*)_k.
$$
We assume that asymptotically the $x_i^*$'s are independent from the $B_{ki}$'s, an assumption supported by the chaos hypothesis, see for instance Geman and Hwang \cite{geman-1982}. Denote by $\mathbb{E}_{\bs{x}^*}=\mathbb{E}(\,\cdot\mid \bs{x}^*)$ the conditional expectation with respect to $\bs{x}^*$. Notice that conditionally to $\bs{x}^*$, the $\check{Z}_k$'s are independent Gaussian random variables, whose two first moments can easily be computed, see Appendix \ref{app:proof-heuristics}, Section \ref{app:proof-heur1} for the details:  
$$
\mathbb{E}_{\bs{x}^*} \check{Z}_k= \mu\,\phat\, \mhat\quad \text{and}\quad \textrm{var}_{\bs{x}^*}(\check{Z}_k)= \frac{\phat\shat^2}{\alpha^2} \, .
$$ 
Notice that the fact that $\mathbb{E}_{\bs{x}^*}$ and $\textrm{var}_{\bs{x}^*}(\check{Z}_k)$ only depend on $\phat, \shat$ and $\mhat$ which are (supposedly) converging quantities supports the idea that $\check{Z}_k$ is unconditionally a Gaussian random variable with moments:
$$
\mathbb{E} \check{Z}_k= \mu\,p^*\, m^*\quad \text{and}\quad \textrm{var}(\check{Z}_k)= \frac{p^* (\sigma^*)^2}{\alpha^2} \, ,
$$
where $p^*,m^*,\sigma^*$ are resp. the limits of $\phat, \mhat, \shat$. We now introduce the standard Gaussian random variables $(Z_k)_{k\in [n]}$ where
$$
Z_k= \frac{\check{Z}_k - \mathbb{E} \check{Z}_k}{\sqrt{\textrm{var}(\check{Z}_k)}} = \alpha \frac{\check{Z}_k - \mu\,p^*\, m^*}{\sigma^*  \sqrt{p^*}}\, .
$$
Consider the equilibrium $\bs{x}^*=(x_k^*)_{i\in [n]}$. If $k\in {\mathcal S}$, that is $x_k^*>0$, we have
$$
1- x_k^* + (B\bs{x}^*)_k=0 \quad \Rightarrow \quad 1+ (B\bs{x}^*)_k>0 \, .
$$ 
This identity has two implications:
$$
x_k^*= 1 + (B\bs{x}^*)_k\quad \textrm{and}\quad 1+ (B\bs{x}^*)_k>0 \qquad \textrm{if}\ k\in {\mathcal S}\, .
$$
Relying on the representation $(B\bs{x}^*)_k=\check{Z}_k$, we obtain the representation 
\begin{equation}\label{eq:representation-x_k}
    x_k = 1+(B\bs{x}^*)_k \ =\ 1+\mu\,p^* \,m^*+\frac{\sigma^*\sqrt{p^*}}{\alpha} Z_k \qquad \textrm{if}\quad k\in {\mathcal S}\,.
\end{equation}
and the condition:
$$
1+(B\bs{x}^*)_k \ =\ 1+\mu\,p^* \,m^*+\frac{\sigma^*\sqrt{p^*}}{\alpha} Z_k\, >\, 0\,.
$$
If $k\notin {\mathcal S}$ then 
$$
1+ (B\bs{x}^*)_k\ =\ 1+\mu\,p^*\,m^* +\frac{\sigma^*\sqrt{p^*}}{\alpha} Z_k\ \le\ 0
$$
by the non invadability condition. Otherwise stated, 
$$
\left\{
\begin{array}{lcll}
 Z_k&\le& -\frac{\alpha(1+\mu\, p^* m^*)}{\sigma^*\sqrt{p^*}} &\text{if}\  k\notin {\mathcal S}\,,\\
Z_k&>& -\frac{\alpha(1+\mu\, p^* m^*)}{\sigma^*\sqrt{p^*}} &\text{if}\ k\in {\mathcal S}\,.\\ 
\end{array}
\right.
$$
Considering the order statistics of the $Z_k$'s we obtain:
$$
Z_1^*\le  \cdots \le Z_{\bs{i}}^*\le -\frac{\alpha(1+\mu\, p^* m^*)}{\sigma^*\sqrt{p^*}} \le Z_{\bs{i}+1}^*\le \cdots \le Z_n^*\, . 
$$
Now, there are exactly $n-|{\mathcal S}|=n(1-\phat)$ indices before the threshold corresponding to the components of $\bs{x}^*$ equal to zero. In particular, index $\bs{i}=n(1-\phat)$ corresponds to the value
$$
Z^*_{\bs{i}}\simeq -\frac{\alpha(1+\mu\,p^*\,m^*)}{\sigma^*\sqrt{p^*}}
$$ 
Relying on \eqref{eq:typical}, we finally obtain
$$
\Phi^{-1}(1-\phat) = -\frac{\alpha(1+\mu\,p^* \, m^*)}{\sigma^*\sqrt{p^*}}\ .
$$
It remains to replace $\phat$ by its limit $p^*$ to obtain \eqref{eq:heur1}.




\paragraph*{Equation \eqref{eq:heur2}.} Our starting point is the following generic representation of an abundance at equilibrium (either of a surviving or vanishing species):
\begin{multline*}
x_k^* = \left( 1+\mu\,p^* m^* +\frac{\sigma^*\sqrt{p^*}}{\alpha} Z_k \right)\bs{1}_{\{Z_k> - \delta^*\}}\\
= \left( 1+\mu\,p^* m^* \right)\bs{1}_{\{Z_k> - \delta^*\}} +\left( \frac{\sigma^*\sqrt{p^*}}{\alpha} Z_k \right)\bs{1}_{\{Z_k> - \delta^*\}}
\end{multline*}
Summing over ${\mathcal S}$ and normalizing, 
\begin{equation*}
\begin{split}
\frac{1}{|{\mathcal S}|}\sum_{k \in \mathcal{S}}x_k^{*} &=(1+\mu\, p^* m^*)\frac{1}{|{\mathcal S}|}\sum_{k \in \mathcal{S}}\bs{1}_{\{Z_k > -\delta^*\}}+\frac{\sigma\sqrt{p^*}}{\alpha}\frac{1}{|{\mathcal S}|}\sum_{k \in \mathcal{S}}Z_k\bs{1}_{\{Z_k > -\delta^*\}}, \\
\mhat &\stackrel{(a)}= (1+\mu\, p^* m^*)+\frac{\sigma\sqrt{p^*}}{\alpha}\frac{n}{|{\mathcal S}|}\frac{1}{n}\sum_{k \in [n]}Z_k\bs{1}_{\{Z_k > -\delta^*\}}, \\
\mhat &\stackrel{(b)}\simeq (1+\mu\,p^* m^*)+\frac{\sigma\sqrt{p^*}}{\alpha} \frac{1}{\mathbb{P}(Z > -\delta^*)}\mathbb{E}(Z\bs{1}_{\{Z > -\delta^*\}}), \\
\mhat &\simeq (1+\mu\, p^* m^*)+\frac{\sigma\sqrt{p^*}}{\alpha}\mathbb{E}(Z \mid Z > -\delta^*).
\end{split}
\end{equation*}
where $(a)$ follows from the fact that $|{\mathcal S}| = \sum_{k \in \mathcal{S}}\bs{1}_{\{Z_k > -\delta^*\}}$ (by definition of ${\mathcal S}$), $(b)$ from the law of large numbers $\frac 1n \sum_{k\in [n]} Z_k \bs{1}_{\{Z_k>-\delta\}} \xrightarrow[n\to\infty]{} \mathbb{E}Z \bs{1}_{\{Z>-\delta\}}$ and $\frac{|{\mathcal S}|}n \xrightarrow[n\to\infty]{} \mathbb{P}(Z>-\delta^*)$ with $Z\sim{\mathcal N}(0,1)$. It remains to replace $\mhat$ by its limit $m^*$ to obtain \eqref{eq:heur2}.

Eq. \eqref{eq:heur3} can be obtained similarly. Details are provided in Appendix \ref{app:proof-heuristics}, see Section \ref{app:proof-heur3}.

\section{Switching between equilibria: changing interaction strength}
\label{sec:ecological}
In the previous sections, the strength $\alpha$ of the interactions was fixed, cf. equation \eqref{eq:LV2}. However, in nature interactions between species are constantly changing due to e.g. abiotic factors such as temperature, which affect the rate at which individuals forage for prey, etc. Our contribution is rooted in the framework of asymptotic dynamics, but many recent ecological studies highlight the importance of taking into account both transient dynamics (out-of-equilibrium abundance values due to frequent perturbations) and shifts between equilibria due to changing environmental conditions \cite{hastings2001,fukami2011,Nolting2016}. In the sequel, we discuss a more general framework.

\paragraph{The model and intuition.} If we restrict ourselves to a specific environment, a possible ecological interpretation of the fluctuation of interaction strength corresponds to the relationship between the size of the habitat and the probability of contact between individuals from two interacting species (e.g. think of a pool of freshwater containing piscivorous fishes and their prey species -- interactions, be them competition or predation, would be potentially more frequent if the volume of water was reduced). In physics, think of particles in motion in a given volume: if the number of particles and the temperature stay constant, reducing the volume should increase the number of interactions between particles.

From a model standpoint, let $\mu$ be fixed, $\alpha=\alpha(t):\mathbb{R}^+\to (\sqrt{2},\infty)$. Consider 
\begin{equation}\label{eq:LV_period}
\frac{dx_k}{dt} = x_k\, \left( 1 -  x_k + (B_t\boldsymbol{x})_k \right)\ ,\quad k\in [n]\ ,
\end{equation}
where matrix $B_t$ admits the following representation
\begin{equation*}
    B_t = \frac{A}{\alpha(t) \sqrt{n}}+\frac{\mu}{n}\bs{1}\bs{1}^T\quad \textrm{and}\quad (\alpha(t),\mu) \in \mathcal{A}\, .
\end{equation*}
\begin{rem} Following Theorem \ref{th:unicite}, condition $(\alpha(t),\mu) \in \mathcal{A}$ guarantees that there exists a unique equilibrium $\bs{x}^*(t)$ for every $t\in \mathbb{R}^+$. 
\label{rem:transition_state}
\end{rem}

\begin{figure}[h!]
\centering
\begin{subfigure}[b]{0.48\textwidth}
  \centering
  \includegraphics[width=\textwidth]{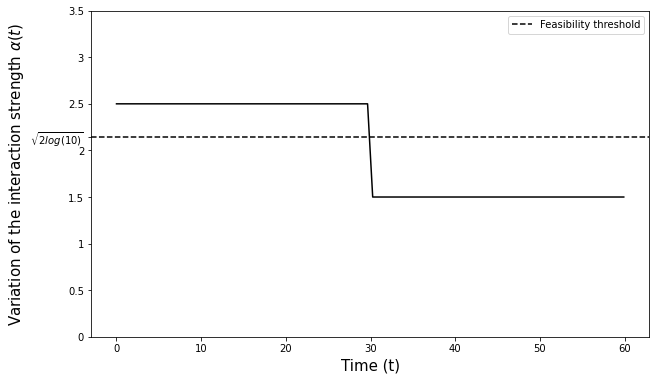}
  \caption{Step function $\alpha(t)$}
  \label{subfig:alpha_var}
\end{subfigure}
\hfill
\begin{subfigure}[b]{0.48\textwidth}
  \centering
  \includegraphics[width=\textwidth]{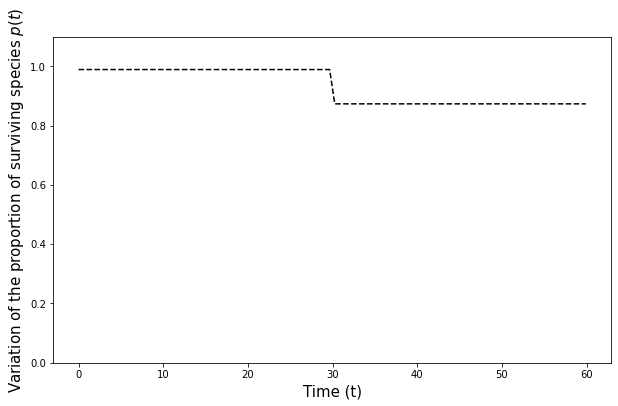}
  \caption{Proportion of surviving species $p(t)$}
  \label{subfig:prop_surv}
\end{subfigure}
   \caption{(a) Variation of the interaction strength through time, used in the dynamics of a ten-species system depicted in Figure \ref{fig:dyn_abrupt} ($\alpha_1 = 2.5$, $\alpha_2 = 1.5$). The dotted line represents the feasibility threshold associated to the system. \newline (b) Variation of the proportion of surviving species depending of the variation of $\alpha(t)$ in (a).}
\label{fig:alpha_type}
\end{figure}

We focus on the case of a system that fluctuates between two equilibrium points (Figure \ref{subfig:alpha_var}). 
We consider a sudden incident, most often irreversible in the short term, which reduces a portion of habitat, e.g. forest fires. The system transits from a feasible state to a state with vanishing species due to the change of the strength of interactions, modelled by the following {\em step function}:
\begin{equation}
    \label{eq:step-function}
\alpha(t) \ = \ \alpha_1\bs{1}_{[0,t_0)}+ \alpha_2\bs{1}_{[t_0,+\infty)}, \quad (\alpha_1,\alpha_2,t_0) \in (\sqrt{2},+\infty)^2\times \mathbb{R}_+
\end{equation}
The change of model parameter occurs at $t_0$ which causes a change in the strength of the interactions going from a value $\alpha_1$ to $\alpha_2$. One may choose $\alpha_2 < \alpha_1$ and the difference (or ratio) between the two values represents the intensity of the incident.

In large dimension, it is possible to characterize this change by its impact on the number of surviving species in the system \eqref{eq:LV_period}. At a given time t, the proportion of surviving species $p=p(t)\in [0,1]$ can be computed by resolving the system in Heuristics \ref{heur:proportion}. This function, associated to the step function $\alpha$ given in \eqref{eq:step-function}, is represented in Figure \ref{subfig:prop_surv}.

\begin{figure}[h!]
\centering
\includegraphics[width=\textwidth]{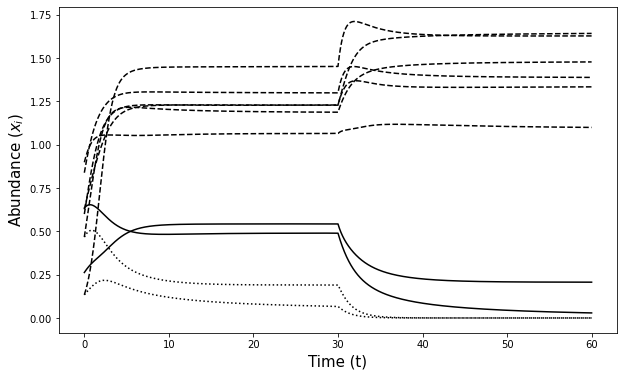}
\caption{Abundance dynamics in the case of a community of ten species. The matrix of interactions $A$  and the initial conditions are common and we apply the function of variation $\alpha(t)$ given in Figure \ref{subfig:alpha_var}. The dashed lines represents species which benefit from habitat variation; solid lines represent species suffering from the change. Dotted lines represent species undergoing extinction.
}
\label{fig:dyn_abrupt}
\end{figure}

In the case of a sudden incident, the proportion of surviving species predicted by the heuristics has a  form similar to $\alpha(t)$ i.e. a step response. In the feasible state, $p(t)$ is close to 1 (i.e. all species coexist); after the transition occurs, some species vanish and here $p(t) \approx 0.87$. Beware that the heuristics results follow instantaneously the change of $\alpha$; however, there is a smoother transition in the dynamics between the two equilibria (respectively corresponding to $\alpha_1$ and $\alpha_2$) due to the return rate to equilibrium, see for instance \cite{neubert_alternatives_1997}, \cite{arnoldi_how_2018}). This transition is illustrated in Figure \ref{fig:dyn_abrupt}.

\paragraph{Simulations.} 
We provide hereafter a simulation-based analysis of the impact of the sudden incident on a given ecosystem: Define a ten-species system \eqref{eq:LV_period} with a fixed matrix of interactions $A$ with Gaussian entries $\mathcal{N}(0,1)$ and consider $\alpha=\alpha(t)$ as in Fig. \ref{subfig:alpha_var}.

This scenario has a mixed impact on the community, see Fig. \ref{fig:dyn_abrupt}. While some species benefit from this change through an increase in their abundances, others are severely affected by this shift, some of which become extinct. This phenomenon can be understood as follows: at first ($t\le t_0$), the system admits a feasible equilibrium state with $\alpha=\alpha_1=2.5>\sqrt{2\log(10)}\simeq 2,14$ and the abundances converge to this equilibrium (see Figure \ref{subfig:alpha_var}).
When the transition occurs at $t=t_0$, Theorem \ref{th:unicite} ensures the convergence to a new equilibrium defined by parameter $\alpha_2$. Since $\alpha_2=1.5$ is below the feasibility threshold $\sqrt{2\log(10)}\simeq 2.14$, some species vanish. In other words, this sudden change of model parameter causes an increase of interaction strengths, which has a negative impact on species diversity. 

\paragraph*{Evolution of diversity} Finally, we illustrate the evolution of diversity using diversity indicators more suited to the description of changes such as the one represented in Fig. \ref{fig:dyn_abrupt} \cite{jost2006entropy}. We introduce here Shannon diversity $H'$, a standard measure of biodiversity in ecology, which is given by 
\begin{equation}
H' = - \sum_{i}\frac{x_i}{\sum_{j}x_j}\log\left(\frac{x_i}{\sum_{j}x_j} \right)
\label{eq:shannon_index}
\end{equation}
and ranges from 0 (one species completely dominates the community) to $\log(n)$, when each species is equally abundant. When many species become rare while others become more abundant, $H'$ drops. Because $H'$ varies before species actually vanish, it is a more sensitive index of community diversity than species richness.

The Hill number of order 1, defined as $e^{H'}$, is a diversity measure that takes into account species abundances and varies between 1 and $n$, i.e. it behaves like an ``effective species richness", see e.g. \cite{jost2007}. 

\begin{figure}[h!]
\centering
\begin{subfigure}[b]{0.48\textwidth}
  \centering
  \includegraphics[width=\textwidth]{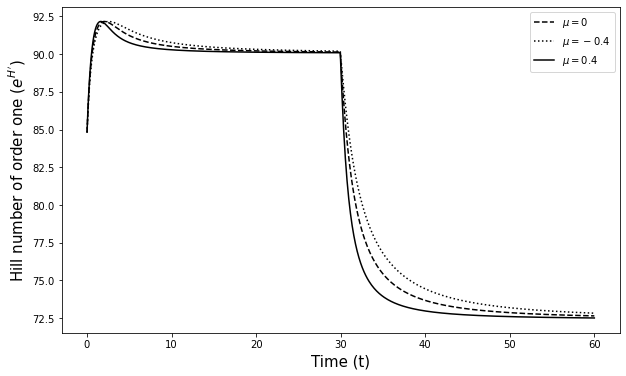}
  \caption{Impact of $\mu$}
  \label{subfig:div_mu}
\end{subfigure}
\hfill
\begin{subfigure}[b]{0.48\textwidth}
  \centering
  \includegraphics[width=\textwidth]{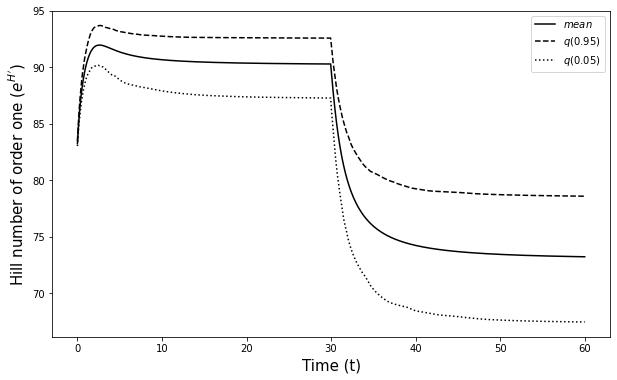}
  \caption{Variability of Hill number $e^{H'}$}
  \label{subfig:div_var}
\end{subfigure}
   \caption{Dynamics of the Hill number of order 1 in the case of an ecosystem of a hundred species. The initial conditions are similar for each species. We define an interaction matrix $A$ and let the dynamics of Lotka-Volterra evolve according to model \eqref{eq:LV_period} and we apply the function of variation $\alpha(t)$ of Figure \ref{subfig:alpha_var}. For each time step, we compute $e^{H'}$. We repeat this scheme a large number of times (here 500), and we average the time series.
   In (a), we apply this procedure for different values of $\mu$. In (b), we apply this procedure for a fixed $\mu = 0$ and compute the quantiles of the 500 trajectories.}
\label{fig:div_index}
\end{figure}

In Fig. \ref{subfig:div_mu}, we represent the mean of this diversity measure over time for a hundred-species system and observe a negative impact of the variation of the strength of the interactions on diversity. Parameter $\mu$ has no impact on diversity at equilibrium (similarly, $\mu$ has no impact on the proportion of persistent species), but the lower the value of $\mu$, the slower the transition to a new equilibrium. In other words, the more generally ``competitive" the ecosystem is (i.e. very negative values of $\mu$), the longer it takes for transient dynamics to settle near equilibria. 

The evolution of the Hill number of order 1 complements the evolution of species richness: when $\alpha$ decreases, the expected number of surviving species decreases (Fig. \ref{fig:alpha_type}); at the same time, $e^{H'}$ decreases even more drastically as the abundance distribution of surviving species gets more heterogeneous. Figure \ref{subfig:div_var} also shows that the variability of the Hill number of order 1 among simulations increases drastically when $\alpha$ decreases. The conclusion is that the more species are lost, the more difficult it is to predict the diversity index as $\sigma^*$ depends on $\alpha$ and strongly influences $e^{H'}$.

\begin{figure}[h!]
\centering
\begin{subfigure}[b]{0.48\textwidth}
  \centering
  \includegraphics[width=\textwidth]{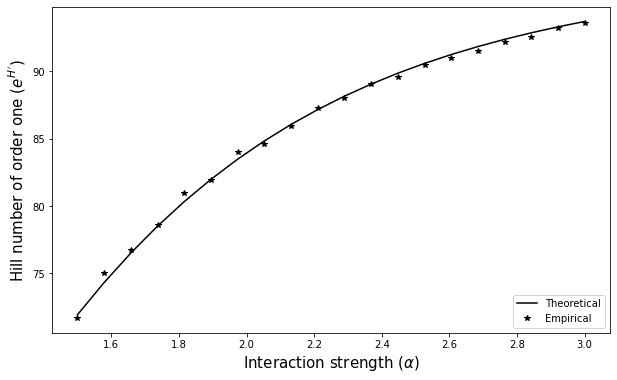}
  \caption{Impact of $\alpha$}
  \label{subfig:TvsE_alpha_div_ind}
\end{subfigure}
\hfill
\begin{subfigure}[b]{0.48\textwidth}
  \centering
  \includegraphics[width=\textwidth]{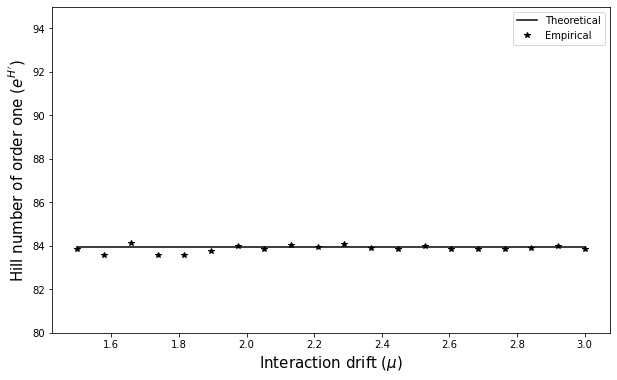}
  \caption{Impact of $\mu$}
  \label{subfig:TvsE_mu_div_ind}
\end{subfigure}
   \caption{Evolution of the Hill number of order 1 as a function of $\alpha$ (a) and $\mu$ (b). The theoretical solutions (solid line) are computed by resolving \eqref{heur:proportion} and integrating the parameter $(p^*,\sigma^*,m^*)$ in equation \eqref{eq:Hill_theo_approx}. The empirical solutions (star marker) are computed by a Monte-Carlo experiment (100 experiments): we define a matrix B of size $100\times 100$, solve the LCP problem and calculate the associated Hill number $e^{H'}$ using \eqref{eq:shannon_index}.}
\label{fig:TvsE_div_ind}
\end{figure}

\paragraph*{Theoretical estimation of diversity}
Standard mathematical methods (Taylor's theorem) can be used to obtain a theoretical approximation of the Hill number of order 1 (details are provided in Appendix \ref{app:div_index}):
\begin{equation}
\label{eq:Hill_theo_approx}
e^{H'} \approx np^*\left(\frac{3}{2}-\frac{1}{2} \frac{(\sigma^{*})^2}{(m^{*})^2}\right)\ .
\end{equation}
This estimator is based on the properties of the persistent species $(p^*,m^*,\sigma^*)$ calculated by solving the fixed point equation of the heuristics \eqref{heur:proportion}. These three properties depend on the type of the interactions, as indicated by parameters $(\alpha,\mu)$ (Figure \ref{fig:test}). We compare the accuracy of this estimator with two examples in which the strength of the interactions $(\alpha)$ and the interaction drift $(\mu)$ vary (Figure \ref{fig:TvsE_div_ind}). 

On the one hand, a shift of the interaction drift $\mu$ does not affect the proportion of surviving species. Furthermore, the impact of $\mu$ on $\sigma^*$ and $m^*$ is proportional i.e. $\frac{\sigma^*}{m^*}$ is equal to a constant. For these reasons, $\mu$ does not affect the Hill number (Figure \ref{subfig:TvsE_mu_div_ind}). On the other hand, if $\alpha$ increases, then $p^*, m^*, \sigma^* \xrightarrow[\alpha \rightarrow +\infty]{} 1$ and $e^{H'} \rightarrow n$ which is makes sense because when $\alpha$ becomes very large, all abundances converge to 1. If $\alpha$ decreases: $p^*$ decreases, and $\sigma^*$ increases faster than $m^*$. This confirms that $e^{H'}$ decreases even more drastically as the abundance distribution of surviving species gets more heterogeneous (Figure \ref{subfig:TvsE_alpha_div_ind}). 

\section{Discussion}

In this paper, our main interest was to describe the impact of the strength $\alpha$ and mean $\mu$ of interactions in large LV models on the conditions of coexistence of interacting species. Combining results from Takeuchi \cite{takeuchi1996global} with standard RMT results, we have provided insights into the study of stability of large random ecosystems - see \cite{stone2018feasibility,gibbs2018effect}. 

We have characterized the unique equilibrium properties of the surviving species by resolving a system of equations. From a physicist point of view, similar equations were obtained by Opper and Diederich \cite{opper_phase_1992} and studied in a more general framework by Bunin using the dynamical cavity method \cite{bunin2017ecological} and Galla \cite{galla_dynamically_2018} using generating functional techniques. 

The coexistence of many species in random ecosystems was also studied by Servan et al. \cite{servan_coexistence_2018} and Pettersson et al. \cite{pettersson_predicting_2020}, where a more generic case was analyzed with different growth rates. Grilli et al. \cite{grilli_feasibility_2017} identified the key quantities regulating the parameter space leading to feasible communities. In contrast to previous approaches, an important feature of our model is the monitoring of the impact of interactions by a normalization factor ($\alpha \sqrt{n}$). From an ecological point of view, one might expect that the larger the number of species, the weaker the interactions will be due to some dilution of interactions among potential interaction partners, which would justify the use of such normalization factors. From a mathematical viewpoint, the normalizing parameter $\alpha$ captures the range of a unique equilibrium and the threshold for feasibility. 


In nature, interactions between species are constantly changing and affected by the environment. Under the assumptions that environmental conditions influence interaction strengths, we have endeavoured in Section \ref{sec:ecological} to study the consequences of a sudden change of environmental conditions, expressed through a decrease in parameter $\alpha$. Solving numerically the Lotka-Volterra system confirms the predictions given by heuristics, i.e. that a decrease in $\alpha$ negatively affects equilibrium species richness. A more precise representation of biodiversity dynamics can be obtained through Hill numbers of order 1 which also decreases after the sudden change in $\alpha$. The dynamics of this diversity measure suggests that the mean of interaction coefficients, $\mu$, affects the duration of transient dynamics, with shorter transient dynamics being associated with more mutualistic interactions (i.e. higher positive values of $\mu$).
\newline

Many questions naturally arise as a follow-up. First, a mathematical proof of the heuristics presented here is a challenging prospect because the LCP procedure induces a statistical dependence that is a priori difficult to handle. However, looking for this proof will certainly help extend the results to other underlying assumptions on the parameters of the LV system.


Regarding the extension of the heuristics to other assumptions, two situations could be of particular interest: non-centered elliptical matrix models as in Bunin \cite{bunin2017ecological} and LV models in which species growth rates are also controlled as in \cite{servan_coexistence_2018}. We are confident that such extensions are possible, but might hinge on more sophisticated developments, in particular to include growth rates in the calculus of order statistics.


In this paper we have only considered the case of a full interaction matrix with parameters $(\alpha,\mu)$. However, food webs are often structured in compartments \cite{baskerville2011} and/or obey hierarchies (e.g. larger species eat smaller ones) \cite{brose2019}. By a numerical analysis of LV systems, one could use the same tools to study more patterned matrices \cite{allesina2012stability}. Recent studies emphasize the sparsity of real food webs \cite{busiello_explorability_2017}. Beyond the feasibility studied by Akjouj and Najim \cite{akjouj2021feasibility}, one could also study the existence and stability of a unique equilibrium in a sparse context.

Finally, variations of the interaction strength highlight the impact of habitat destruction. Many theoretical studies provide mathematical formulas for the return rate to equilibrium \cite{neubert_alternatives_1997,arnoldi_how_2018}. A further theoretical study of model \eqref{eq:LV_period} could provide a more quantitative answer. In this article, we have limited the analysis to the case of a single sudden incident, but other types of fluctuations for the interaction strength could be considered for a better understanding of habitat conservation phenomena. For example, a seasonal model could be appropriate to describe the evolution of the dynamics over the seasons.  

\newpage

\vfill\pagebreak

\bibliographystyle{alpha}
\bibliography{mathematics}

\appendix

\section{Simulation details}
Simulations were performed in Python. All the figures and the code are available on Github \cite{code}.

Simulations on the properties of surviving species are performed in two different ways. The theoretical solutions are obtained resolving numerically the system of equations of heuristics \ref{heur:proportion}. We use a solver (cf. scipy.optimize) to find a local minimum of the function defined by the system of equations (a modification of the Powell hybrid method). The empirical solutions are computed using a Monte Carlo experiment. We simulate a large number of matrix matrix B, we resolve the associated LCP problem using the Lemke's algorithm. Then, we use the LCP solution to calculate the properties of the surviving species: proportion of survivors, etc. Finally, we make an average on the ensemble of experiments. The Lemke algorithm is implemented in the lemkelcp package and can be found on Github \cite{packagelemke}. The dynamics of the Lotka-Volterra are achieved by a Runge-Kutta of order 4 (RK4) implemented in the code.

\section{Remaining computations}
\label{app:proof-heuristics}

\subsection{Details on Equation \eqref{eq:heur1}: Moments of $\check{Z}_k$.}
\label{app:proof-heur1}
We compute hereafter the conditional mean and variance of $\check{Z}_k=(B\bs{x}^*)_k$ with respect to $\bs{x}^*$.
We rely on the following identities:
$$
\mathbb{E} B_{ki} =\frac \mu n\,,\quad \mathbb{E} (B_{ki})^2 = \frac 1{\alpha^2 n } + \frac{\mu^2}{n^2} \simeq \frac 1{\alpha^2 n}\,,\quad \mathbb{E} B_{ki}B_{kj} = \frac{\mu^2}{n^2}\quad (i\neq j)\, .
$$
We first compute the conditional mean:
$$
    \mathbb{E}_{\bs{x}^*}(\check{Z}_k) = \sum_{i\in[n]}  \mathbb{E}(B_{ki}) x_i^*
    = \sum_{i\in {\mathcal S}}  \mathbb{E}(B_{ki}) x_i^* =  \frac{\mu}{n}\sum_{i\in {\mathcal S}} x_i^*, \\
    = \mu \frac{|{\mathcal S}|}n \frac{1}{|{\mathcal S}|}\sum_{i\in {\mathcal S}} x_i^*, \\
    = \mu\, \phat\, \mhat\, .
$$
We now compute the second moment:
\begin{eqnarray*}
    \mathbb{E}_{\bs{x}^*}(\check{Z}_k^2) &=&\mathbb{E}_{\bs{x}^*} \left( \sum_{i\in {\mathcal S}}  B_{ki} x_i^{*}\right)^2
    \ =\  \mathbb{E}_{\bs{x}^*}  \sum_{i,j\in {\mathcal S}} B_{ki}B_{kj}x_i^*x_j^*\,,\\ 
    &=& \sum_{i \in \mathcal{S}}\mathbb{E}(B_{ki}^2)x_i^{*2}+\sum_{i \neq j}\mathbb{E} (B_{ki}B_{kj})x_i^*x_j^*, \\
    &= &\frac{1}{\alpha^2 n} \sum_{i \in \mathcal{S}}x_i^{*2}+\frac{\mu^2}{n^2} \sum_{i \neq j}x_i^*x_j^*,\\
    &\stackrel{(a)}\simeq& \frac{\phat \shat^2}{\alpha^2}+\frac{\mu^2\, \phat^2}{|{\mathcal S}|^2}\sum_{i,j\in {\mathcal S}}x_i^*x_j^* 
    \ =\ 
        \frac{\phat \shat^2}{\alpha^2}+\frac{\mu^2\, \phat^2}{|{\mathcal S}|^2}\left(\sum_{i\in {\mathcal S}}x_i^*\right)^2
        \ =\ \frac{\phat\shat^2}{\alpha^2}+\mu^2\, \phat^2 \mhat^2\, ,
\end{eqnarray*}
where the approximation in $(a)$ follows from the fact that 
$$
\frac 1{|{\mathcal S}|^2} \sum_{i,j\in {\mathcal S}}x_i^*x_j^* = \frac 1{|{\mathcal S}|^2} \sum_{i \neq j}x_i^*x_j^* + {\mathcal O} \left( \frac 1{|{\mathcal S}|} \right)\, .
$$
We can now compute the variance:
$$
\textrm{var}_{\bs{x}^*} \left(\check{Z}_k\right) =   \mathbb{E}_{\bs{x}^*}\left(\check{Z}_k^2\right) -  \left( \mathbb{E}_{\bs{x}^*}\check{Z}_k\right)^2 \ =\ \frac{\phat \, \shat^2}{\alpha^2}\, .
$$

\subsection{Details on Equation \eqref{eq:heur3}.}
\label{app:proof-heur3}
As for the proof of \eqref{eq:heur2}, we start from the generic representation of $x_k^*$:
\begin{multline*}
x_k^*\ =\ \left(1+\mu\, p^* m^* +\frac{\sigma^*\sqrt{p^*}}{\alpha} Z_k\right)\bs{1}_{\{Z_k > -\delta^*\}}\\
=\ \left(1+\mu\,p^* m^*\right)\bs{1}_{\{Z_k > -\delta^*\}}+\frac{\sigma\sqrt{p^*}}{\alpha}Z_k\bs{1}_{\{Z_k > -\delta^*\}}\, .
\end{multline*}
Taking the square, we get 
\begin{multline*}
x_k^{*2} = \left(1+\mu\,p^* m^*\right)^2\bs{1}_{\{Z_k > -\delta\}} \\+2(1+\mu\,p^* m^*)\frac{\sigma^*\sqrt{p^*}}{\alpha}Z_k\bs{1}_{\{Z_k > -\delta\}}
+\frac{(\sigma^*)^2p^*}{\alpha^2}Z^2_k\bs{1}_{\{Z_k > -\delta^*\}}\, .
\end{multline*}
Summing over ${\mathcal S}$ and normalizing, we get 
\begin{multline*}
    \frac{1}{|{\mathcal S}|}\sum_{k \in \mathcal{S}}(x_k^*)^{2} = (1+\mu\, p^* m^*)^2\frac{1}{|{\mathcal S}|}\sum_{k \in \mathcal{S}}\bs{1}_{\{Z_k > -\delta^*\}} \\
    +2(1+\mu\,p^* m^*)\frac{\sigma^*\sqrt{p^*}}{\alpha}\frac{1}{|{\mathcal S}|}\sum_{k \in \mathcal{S}}Z_k\bs{1}_{\{Z_k > -\delta^*\}}\\+\frac{(\sigma^*)^2 p^*}{\alpha^2}\frac{1}{|{\mathcal S}|}\sum_{k \in \mathcal{S}}Z^2_k\bs{1}_{\{Z_k > -\delta^*\}}\,.
\end{multline*}
Finally, we conclude by replacing the empirical means by their limits
\begin{eqnarray*}
\frac{1}{|{\mathcal S}|}\sum_{k \in \mathcal{S}}Z^i_k\bs{1}_{\{Z_k > -\delta^*\}} &=& \mathbb{E}(Z^i\mid Z>-\delta^*)\ ,\quad i=1,2\, .
\end{eqnarray*}
and get 
\begin{multline*}
\shat^2 = (1+\mu\,p^* m^*)^2+2(1+\mu\,p^* m^*)\frac{\sigma^*\sqrt{p^*}}{\alpha}\mathbb{E}(Z \mid Z > -\delta^*)\\+\frac{(\sigma^*)^2p^*}{\alpha^2}\mathbb{E}(Z^2 \mid Z > -\delta^*)\,.
\end{multline*}
It remains to replace $\shat$ by its limit $\sigma^*$ to obtain \eqref{eq:heur3}.

\subsection{Density of the distribution of the persistent species.}
Assume that $x^*>0$, and les $f=\mathbb{R}\to \mathbb{R}$ be a bounded continuous test function. We have
\begin{eqnarray*}
    \mathbb{E}f(x_k^*) &=& \mathbf{E}\left[f\left(1+ \frac{\sigma^*\sqrt{p^*}}{\alpha}Z_k+\mu\,p^* m^*\right)\  \bigg|\  Z_k > -\delta^* \right]\ , \\
    &=& \int_{-\infty}^{\infty}f\left(1+\mu\,p^* m^*+\frac{\sigma^*\sqrt{p^*}}{\alpha}u\right) \frac{\bs{1}_{\{u>-\delta^*\}}}{1-\Phi(-\delta^*)}\frac{e^{-\frac{u^2}{2}}}{\sqrt{2p^*}}du\ , \\
    &=& \int_{0}^{\infty}f(y) e^{-\frac 12\left(\frac{\alpha}{\sigma^* \sqrt{p^*}}y-\delta^*\right)^2}\frac{\alpha}{\sqrt{2} \Phi(\delta^*)\,p^*\,  \sigma^*}\,dy\ ,
\end{eqnarray*}
hence the density of $x^*_k$.

\subsection{Theoretical estimation of the diversity index}
\label{app:div_index}
Recall that $|\mathcal{S}| = n \hat{p}$ is the number of surviving species and that 
$$p_i = \frac{x_i}{\sum_{j\in {\mathcal S}} x_j}$$
is the frequency of (surviving) species $i$. 

To find a theoretical estimate of Hill number of order 1, we proceed by expansion and set $$p_i = \frac{1}{|\mathcal{S}|}+\delta_i\,,\quad |\delta_i|\ll \frac 1{|S|}$$ where $\delta_i$ represents the deviation of species $i$ from the standard frequency if all surviving species have the same frequency. Notice that $\sum_{i\in \mathcal{S}}\delta_i=0$.
$$
    H' \ =\ - \sum_{i\in \mathcal{S}} p_i \log(p_i) 
    = - \sum_{i\in \mathcal{S}} \left(\frac{1}{|\mathcal{S}|}+\delta_i \right)\log\left(\frac{1}{|\mathcal{S}|}+\delta_i \right) \ .
$$
We use the Taylor-Young formula of order 2 to decompose the log: 
\begin{eqnarray*}
  \log\left(\frac{1}{|\mathcal{S}|}+\delta_i \right) &=& \log \left(\frac{1}{|\mathcal{S}|}\right)+|\mathcal{S}|\delta_i-\frac{\delta_i^2|\mathcal{S}|^2}{2}+\delta_i^3\varepsilon(\delta_i) \ , \\
  &\approx& \log \left(\frac{1}{|\mathcal{S}|}\right)+|\mathcal{S}|\delta_i-\frac{\delta_i^2|\mathcal{S}|^2}{2} \ .
\end{eqnarray*}
\begin{eqnarray*}
   H' &\approx&  - \sum_{i\in \mathcal{S}} \left(\frac{1}{|\mathcal{S}|}+\delta_i \right) \left(\log \left(\frac{1}{|\mathcal{S}|}\right)+|\mathcal{S}|\delta_i-\frac{\delta_i^2|\mathcal{S}|^2}{2}\right)\ , \\
   &=& - \sum_{i\in \mathcal{S}}\left[\frac{1}{|\mathcal{S}|}\log \left( \frac{1}{|\mathcal{S}|}\right)+\delta_i-\frac{\delta_i^2|\mathcal{S}|}{2}+\delta_i \log \left(\frac{1}{|\mathcal{S}|} \right)+ |\mathcal{S}|\delta_i^2-\frac{\delta_i^3 |\mathcal{S}|^2}{2} \right]\ , \\
   &=& \log(|\mathcal{S}|) - \sum_{i\in \mathcal{S}} \frac{\delta_i^2|\mathcal{S}|}{2} +\sum_{i\in \mathcal{S}} \frac{\delta_i^3|\mathcal{S}|^2}{2}\ .
\end{eqnarray*}
Notice that $\sum_{i=1}^{|\mathcal{S}|} \frac{\delta_i^3|\mathcal{S}|^2}{2}$ is negligible since $|\delta_i|\ll |S|^{-1}$.
The term 1 corresponds to the maximum value that the Shannon diversity index can take if $|\mathcal{S}|$ are present in the system. It remains to develop the second term of the r.h.s.
\begin{eqnarray*}
-\frac{1}{2}\sum_{i\in \mathcal{S}}\delta_i^2|\mathcal{S}| 
&=& -\frac{|\mathcal{S}|}{2}\sum_{i\in \mathcal{S}} \left(\frac{x_i}{\sum_{j\in \mathcal{S}} x_j}-\frac{1}{|\mathcal{S}|} \right)^2\ , \\
&=& -\frac{|\mathcal{S}|}{2} \sum_{i\in \mathcal{S}} \left(\frac{x_i^2}{(\sum_{j\in \mathcal{S}}x_j)^2}-\frac{2}{|\mathcal{S}|}\frac{x_i}{\sum_{j\in \mathcal{S}} x_j}+\frac{1}{|\mathcal{S}|^2} \right)\ , \\
&=& -\frac{|\mathcal{S}|}{2} \sum_{i\in \mathcal{S}} \left(\frac{x_i^2}{(\sum_{j\in \mathcal{S}}x_j)^2}\right) +\frac{1}{2}\ , \\
&=& -\frac{|\mathcal{S}|}{2}\frac{\sum_{i\in \mathcal{S}}x_i^2}{|\mathcal{S}|^2 (\frac{1}{|\mathcal{S}|}\sum_{j\in \mathcal{S}}x_j)^2}+\frac{1}{2}\ ,\\
&=& -\frac{1}{2}\frac{\frac{1}{|\mathcal{S}|}\sum_{i\in \mathcal{S}}x_i^2}{ (\frac{1}{|\mathcal{S}|}\sum_{j\in \mathcal{S}}x_j)^2}+\frac{1}{2}\ , \\
&=&-\frac{1}{2} \frac{\hat{\sigma}^2}{(\hat{m})^2}+\frac{1}{2}\ , \\
&=& -\frac{1}{2}\left(\frac{\hat{\sigma}^2}{\hat{m}^2}-1  \right) \ .
\end{eqnarray*}
Finally the Hill number of order 1 can be computed as:
\begin{eqnarray*}
    e^{H'} &\approx& e^{\log(|\mathcal{S}|)-\frac{|\mathcal{S}|}{2}\sum_{i=1}^{|\mathcal{S}|}\delta_i^2}\ , \\
     &\approx& |\mathcal{S}|\left(1-\frac{|\mathcal{S}|}{2}\sum_{i=1}^{|\mathcal{S}|}\delta_i^2\right)\ =\
     |\mathcal{S}|\left(1-\frac{1}{2} \frac{\hat{\sigma}^2}{(\hat{m})^2}+\frac{1}{2}\right)\ =\
      \frac{|\mathcal{S}|}2\left(3- \frac{\hat{\sigma}^2}{(\hat{m})^2}\right)\ . 
\end{eqnarray*}
Replacing $|{\mathcal S}|$ by $np^*$ and $\hat{\sigma}$ and $\hat{m}$ by their limits, we get the desired result:
\begin{eqnarray*}
    e^{H'} &\approx&      \frac{np^*}2\left(3- \frac{(\sigma^*)^2}{(m^*)^2}\right)\ . 
\end{eqnarray*}

\end{document}